\documentclass[12pt]{article}
\usepackage{amsfonts}
\usepackage{latexsym, amssymb, amsmath, amscd, amsfonts, epsfig, graphicx}
\usepackage{mathrsfs}
\usepackage{colordvi}
\usepackage{ifpdf}
\newtheorem{theorem}{Theorem}[section]

\newtheorem{lemma}[theorem]{Lemma}

\newtheorem{corollary}[theorem]{Corollary}

 \numberwithin{equation}{section}

\def\qed{\nopagebreak\hfill{\rule{4pt}{7pt}}}
\def\proof{\noindent {\it{Proof.} \hskip 2pt}}

\def\pf{\noindent {\it Proof.} }

\begin{document}

\begin{center}
{\large\bf Minimal Permutations and $2$-Regular Skew Tableaux}
\end{center}

\begin{center}
William Y.C. Chen$^{1}$, Cindy C.Y. Gu$^{2}$ and Kevin J. Ma$^{3}$ \\[6pt]
Center for Combinatorics, LPMC-TJKLC\\
Nankai University, Tianjin 300071, P. R. China

Email: $^{1}${\tt chen@nankai.edu.cn}, $^{2}${\tt
guchunyan@cfc.nankai.edu.cn}, \\$^{3}${\tt majun@cfc nankai.edu.cn},

\end{center}

\vspace{0.3cm} \noindent{\bf Abstract.}
Bouvel and Pergola introduced the notion of
 minimal permutations in the study of
 the whole genome
  duplication-random loss model for genome rearrangements.
 Let
  $\mathcal{F}_d(n)$ denote the set of minimal
 permutations of length $n$ with $d$ descents, and let
 $f_d(n)= |\mathcal{F}_d(n)|$.
They derived that
 $f_{n-2}(n)=2^{n}-(n-1)n-2$ and $f_n(2n)=C_n$,
 where $C_n$ is the $n$-th Catalan number.
 Mansour and Yan proved that
  $f_{n+1}(2n+1)=2^{n-2}nC_{n+1}$.
  In this paper, we consider the problem of
  counting minimal permutations in $\mathcal{F}_d(n)$ with
  a prescribed set of ascents. We show that such structures
   are in one-to-one correspondence with a class of
   skew Young tableaux, which we call $2$-regular
   skew tableaux. Using the determinantal formula for the number
    of skew Young tableaux of a given shape, we find an explicit
    formula for $f_{n-3}(n)$.
     Furthermore, by using the Knuth equivalence, we give a combinatorial interpretation of a formula for a refinement
     of the number $f_{n+1}(2n+1)$.

\vskip 3mm

\noindent {\bf Keywords:} minimal permutation,
 $2$-regular skew tableau, Knuth equivalence,
 the RSK algorithm.

\vskip 3mm

\noindent {\bf AMS Classification:} 05A05, 05A19.

\section{Introduction}

The notion of minimal permutations was introduced by Bouvel and Pergola
in the study of genome evolution, see \cite{bou}.
Such permutations are a basis of permutations that can be obtained from the identity permutation
via a given number of steps in the duplication-random loss model.
Let $\pi=\pi_1\pi_2\cdots \pi_n$ be
a permutation. A duplication of $\pi$ means the duplication
 of a fragment of consecutive elements of $\pi$ in such a way that the duplicated
fragment is put immediately after the original fragment.
Suppose that $\pi_i\pi_{i+1}\cdots \pi_{j}$ is the fragment
for duplication, then the duplicated sequence is
\[ \pi_1\cdots\pi_{i-1} \pi_i \cdots \pi_j \pi_i \cdots \pi_{j} \pi_{j+1}\cdots \pi_n.\]
A random loss means to
randomly delete one occurrence of each repeated element
$\pi_k$ for $i\leq k \leq j$, so that
we get a permutation again. In the following example,
the fragment $234$ is duplicated, and the underlined elements
are the occurrences of repeated elements that are supposed to
be deleted,
$$1\overbrace{234}56\rightsquigarrow 1\overbrace{234}\overbrace{234}56\rightsquigarrow 1\underline{2}3\underline{4}2\underline{3}456\rightsquigarrow 132456.$$

To describe the notation of minimal permutations, we give an
overview of the descent set of a permutation and the patterns
of subsequences of a permutation.
Let $S_n$ be the set of permutations on $[n]=\{1,2,\ldots,n\}$, where $n\geq 1$. In a permutation $\pi=\pi_1\pi_2\cdots \pi_n\in S_n$, a descent is a position $i$ such that $i\leq n-1$ and $\pi_i>\pi_{i+1}$, whereas an ascent is a position $i$ with $i\leq n-1$ and $\pi_i<\pi_{i+1}$. For example, the permutation $3145726\in S_7$ has two descents $1$ and $5$ and has four ascents  $2,3,4$ and $6$.

Let $V=\{v_1,v_2,\ldots,v_n\}$
be a set of distinct integers listed in increasing order,
namely, $v_1<v_2<\cdots <v_n$. The standardization of a permutation $\pi$ on $V$ is the permutation $\rm{st}(\pi)$ on $[n]$ obtained from $\pi$ by replacing  $v_i$ with  $i$. For example, $\rm{st}(9425)=4213$. A subsequence $\omega=\pi_{i(1)}\pi_{i(2)}\cdots \pi_{i(k)}$ of $\pi$ is said to be of  type $\sigma$ or $\pi$ contains a pattern $\sigma$ if $\rm{st}(\omega)=\sigma$.
 We say that a permutation $\pi\in S_n$ contains a pattern
 $\tau\in S_k$ if there is a subsequence of $\pi$ that is
  of type $\tau$. For example, let $\pi=263751498$.
  The subsequence $3549$ is of type $1324$, and so $\pi$
   contains the pattern $1324$. We use the notation
   $\tau\prec \pi$ to denote that a permutation
   $\pi$ contains the pattern $\tau$, and we use $S_n(\tau_1,
    \ldots, \tau_k)$ to denote the set of permutations $\pi\in S_n$ that avoid the patterns $\tau_1, \tau_2,\ldots, \tau_k$.

A permutation $\pi$ is called a minimal permutation with $d$ descents if it is minimal in the sense that
there exists no permutation $\sigma$ with exactly $d$ descents such that $\sigma\prec\pi$. Denote by $\mathcal {B}_d$ the set of minimal permutations with $d$ descents. Bouvel and Pergola \cite{bou} have shown that the length, namely, the number of elements,
 of any minimal permutation in the set $\mathcal {B}_d$ is at least $d+1$ and at most $2d$. They also proved that in the
whole genome duplication-random loss model,
the permutations that can be obtained from the identity permutation in at most $p$ steps can be characterized as
permutations  $d = 2^p$ descents that avoid certain patterns.

\begin{theorem}[Bouvel and Pergola] \label{con}
Let $\pi=\pi_1\pi_2\cdots\pi_n$ be a permutation on $[n]$. Then $\pi$ is a minimal permutation with $d$ descents if and only if $\pi$ is a permutation with $d$ descents satisfying the following
conditions:
\begin{itemize}
\item[(1)]  It starts and ends with a descent;
\item[(2)] If $i$ is an ascent, that is, $\pi_i<\pi_{i+1}$,
 then $i\in\{2,3,\ldots,n-2\}$ and $\pi_{i-1}\pi_i\pi_{i+1}\pi_{i+2}$ is of type $2143$ or $3142$.
\end{itemize}
\end{theorem}
 Denote by $\mathcal {F}_d(n)$ the set of minimal permutations of length $n$ with $d$ descents and $f_d(n)=|\mathcal{F}_d(n)|$. Clearly, $f_d(n)=0$ for all $d\leq 0$ or $d \geq n$, and $f_d(d+1)=1$ for all $d\geq 1$.
  Bouvel and Pergola  proved that $f_n(2n)$ equals the $n$-th Catalan number, that is,
$$f_n(2n)=C_n=\frac{1}{n+1}{2n\choose n}$$
and $f_{n-2}(n)$  is given by the formula
$$f_{n-2}(n)=2^{n}-(n-1)n-2.$$
  Mansour and Yan  \cite{man} have shown that
\begin{equation}\label{mn}
f_{n+1}(2n+1)=2^{n-2}nC_{n+1}.
\end{equation}

As mentioned by Bouvel and Pergola that it is an open problem
to compute  $f_d(n)$ for other cases of $d$.
In this paper, we consider the enumeration
 of minimal permutations in $\mathcal{F}_d(n)$ with
  a prescribed set of ascents. We show that such minimal
  permutations  are in one-to-one correspondence with a class of
   skew Young tableaux, which we call $2$-regular
   skew tableaux. As a result, we may
    employ the determinant  formula for the number
     of skew Young tableaux of a given shape to
     compute the number $f_d(n)$. With this method, we can unite the known results.
 Moreover, we derive an explicit
    formula for $f_{n-3}(n)$.

   For the number $f_{n+1}(2n+1)$, we obtain a refined formula
   from the determinant formula. Moreover, we give a
   combinatorial interpretation of this  formula by using the
   Knuth equivalence of permutations.

\section{$2$-Regular skew tableaux}

In this section, we establish a connection
between the minimal permutations and skew Young tableaux
of certain shape. To describe our correspondence,
let us give an overview of necessary terminology on Young tableaux as used in Stanley \cite{sta}.

 A partition of a positive integer $n$ is defined to be a sequence $\lambda=(\lambda_1,\ldots,\lambda_k)$ of positive integers
 such that $\sum\lambda_i=n$ and $\lambda_1\geq \cdots \geq\lambda_k$.  If $\lambda$ is a partition of $n$, we write $\lambda\vdash n$, or $|\lambda|=n$. The Ferrers diagram
 of a partition $\lambda$ is a diagram with left-justified rows in which
 the $i$-th row consists of $\lambda_i$ dots. The conjugate partition $\lambda'$ of
 $\lambda$ is obtained by transposing the Ferrers diagram of $\lambda$. The positive terms $\lambda_i$ are called the parts of $\lambda$, and the number of parts is denoted by $l(\lambda)$.

A standard Young tableau (SYT) on $[n]$ is said to be
 of size $n$. If $\lambda$ and $\mu$ are partitions with $\mu\subseteq \lambda$, namely, $\mu_i\leq \lambda_i$
 for all $i$, we can define a standard tableau of skew shape $\lambda/\mu$ as a tableau on $[n]$ that is
 increasing in every row and every column.  The number of
 boxes of the Young diagram of shape $\lambda/\mu$ is denoted
 by $|\lambda/\mu|$.
  For example, below are an SYT  of shape $(4,3,3,1)$ and a skew Young tableau of shape $(6,5,2,2)/(3,1)$:
\[
\begin{array}{ccccccccccccccc}
1 & 3 & 5 & 6 &&&&&&&&& 7& 8& 11   \\
2 & 4 & 8 &&&&&&&& 1& 5 & 9 & 10\\
7 & 9 & 11  &&&&&&& 2&4\\
10  &&&,&&&&&& 3&6
\end{array}.
\]
Recall that if $|\lambda/\mu|=n$ and $l(\lambda)= r$, then the number  of skew Young tableaux of shape $\lambda/\mu$ is given by
\begin{equation}\label{formula}
f^{\lambda/\mu}=n!\,\det\left(\frac{1}
{(\lambda_i-\mu_j-i+j)!}\right)_{i,j=1}^r,
\end{equation}
see, for example, \cite[Corollary 7.16.3]{sta}.

 Let $\{a_1,\ldots,a_k\}$ be a sequence of positive
 integers such that $a_i\geq 2$ for all $i$ and $a_1+a_2+\cdots +a_k=n$.
 Let $P$ be a skew Young tableau of size $n$ with column lengths $a_1,a_2,\ldots,a_k$. We say that
 $P$ is $2$-regular if any two consecutive columns
 overlap exactly by two rows, namely, for any two consecutive columns there are exactly two rows containing elements in
 both columns.
Denote by $\mathcal{P}_{a_1,a_2,\ldots,a_k}(n)$
  the set of $2$-regular skew tableaux with column lengths $a_1,a_2,\ldots,a_k$.

For example, the following skew Young tableau is
$2$-regular and it belongs to  $\mathcal{P}_{4,2,5,3,2}(16)$:
\begin{equation}\label{ol}
\begin{array}{cccccccc}
&&&6&8\\
&&2  & 10 & 15\\
&&5 & 11\\
&& 9\\
1&3 &12\\
4&7&14\\
13\\
16
\end{array}.
\end{equation}

For a permutation $\pi$ of length $n$, a substring of $\pi$ is a sequence of
consecutive elements of $\pi$.
A maximal decreasing substring of $\pi$
is defined to be a decreasing substring
that is not a substring of another decreasing substring. For example, the permutation
$5\,2\,7\,3\,1\,4\,8\,9\,6$ contains five maximal decreasing substrings, namely,
$5\,2, 7\,3\,1,4,8$ and $9\,6$.

It is clear that any permutation $\pi$ with $k-1$ ascents
can be decomposed into $k$ maximal decreasing substrings.
To describe the ascent set, we find it convenient to
use a sequence  $(a_1,a_2,\ldots,a_k)$ to denote
 the lengths of the maximal decreasing substrings,
 and this sequence is called the ascent sequence of $\pi$.
Then the ascent set $\pi$ is expressed as $\{a_1,a_1+a_2,\ldots,a_1+a_2+\cdots+a_{k-1}\}$.

\begin{lemma}
Given a minimal permutation $\pi=\pi_1\pi_2\cdots\pi_n$. Suppose $(a_1,a_2,\ldots,a_k)$ is its ascent sequence, then
 $a_i\geq 2$ for all $i$.
\end{lemma}

\pf By condition (i) of Theorem \ref{con}, $\pi$ starts and ends with a descent, this implies that $a_1\geq 2$ and $a_k\geq 2$.
For each ascent $j=a_1+\cdots+a_i$ of $\pi$, where $1\leq i\leq k-1$, the condition (ii) of Theorem \ref{con}  says that  $\pi_{j-1}\pi_{j}\pi_{j+1}
\pi_{j+2}$ is of type $2143$ or $3142$, which means that both $j-1$ and $j+1$ are descents. Therefore, $\pi$ contains no consecutive ascents, and the length of decreasing sequences containing $\pi_{j-1}\pi_{j}$ and $\pi_{j+1}
\pi_{j+2}$ are least  two. This completes the proof. \qed

Let
$\mathscr{F}_{a_1,a_2,\ldots,a_k}(n)$ denote the set of minimal
permutations of length $n$ with the ascent sequence $(a_1,a_2,\ldots,a_k)$,
and let $F_{a_1,a_2,\ldots,a_k}(n)=|\mathscr{F}_{a_1,a_2,\ldots,a_k}(n)|$.
The following theorem asserts that the number of minimal permutations with a prescribed ascent sequence is equal to the number of skew Young tableaux with fixed column lengths.

\begin{theorem}\label{f}
Let $(a_1,a_2,\ldots,a_k)$ be a sequence of positive integers such that $a_1+a_2\cdots +a_k=n$  and $a_i\geq 2$ for all $i$.
Then there exists a bijection between the set
$\mathscr{F}_{a_1,a_2,\ldots,a_k}(n)$ of minimal permutations
with ascent sequence $(a_1,a_2,\ldots,a_k)$ and the set
$\mathcal{P}_{a_1,a_2,\ldots,a_k}(n)$
of $2$-regular skew tableaux with column
lengths $a_1,a_2,\ldots,a_k$.
\end{theorem}

\pf Suppose $\pi=\pi_1\pi_2\cdots\pi_n$ is a minimal permutation in
$\mathscr{F}_{a_1,a_2,\ldots,a_k}(n)$.
Let $p_i=\pi_{a_i+1}\pi_{a_i+2}\cdots \pi_{a_{i+1}}$
($0\leq i\leq k-1$ and set $a_0=0$) be the $k$ maximal
decreasing substrings of $\pi$, then the elements in each
$p_i$ are strictly decreasing.
Furthermore, by Theorem \ref{con}, if $j=a_1+\cdots+a_i$
 is an ascent then $\pi_{j-1}\pi_{j}\pi_{j+1}
\pi_{j+2}$ is of type $2143$ or $3142$. Therefore, if we place these
four elements into an array as follows,
 \begin{equation}\label{ii}
\begin{array}{cc}
\pi_{j} & \pi_{j+2}\\
\pi_{j-1} & \pi_{j+1},
\end{array}
\end{equation}
then both its rows and columns are strictly increasing.
We next construct a tableau $P$ corresponding to $\pi$ as follows.
Place the elements of each maximal decreasing substring
$p_i$ in one single column, with the decreasing order
from the bottom upward. This guarantees that each column
of the tableau $P$ is strictly increasing.
Now, for every two adjacent maximal decreasing
substrings $p_i$ and $p_{i+1}$,
we assume that the last two elements in $p_i$
and the first two elements in $p_{i+1}$ are
arranged into a $2\times 2$ square as exhibited
 in \eqref{ii}.
 This ensures that each row of $P$ is also
 strictly increasing.
Therefore, $P$ is indeed a $2$-regular skew tableau.
To be more precise, $P$ has the following form,
\begin{equation}\label{aa}
P=
\begin{array}{cccccccc}
            &                 &   \vdots  \\
            &                 &   \vdots      & \cdots & \cdots   \\
            &  \pi_{a_1+a_2}        &   \pi_{a_1+a_2+2}  \\
            &  \pi_{a_1+a_2-1}      &   \pi_{a_1+a_2+1}  \\
         &  \vdots  \\
    \pi_{a_1} &\pi_{a_1+2}\\
     \pi_{a_1-1}&\pi_{a_1+1}\\
    \vdots\\
    \pi_2 \\
    \pi_1\\
\end{array}
\end{equation}
For example,
$\pi=\underline{16}\,\underline{13}\,4\,1\,7\,3\,\underline{14}\,\underline{12}
\,9\,5\,2\,\underline{11}\,\underline{10}\,6\,\underline{15}\,8
\in\mathcal{F}_{11}(16)$
contains $5$ maximal decreasing substrings.
The $2$-regular skew tableau corresponding
to $\pi$ is given by the array in \eqref{ol}.

Conversely, given a $2$-regular skew tableau $P$, we write down
the elements of $P$ from bottom up  and from left to
right. Then we obtain a minimal permutation $\pi$.
Thus we get a bijection.
\qed

For example, all the minimal permutations in
$\mathcal{F}_{n}(2n)$ have
alternating isolated descents as well as
 alternating isolated ascents. Note that
 these minimal permutations start and end with
 descents, see \cite{bou}.
 Therefore, the corresponding $2$-regular skew
tableaux always have straight shape $(n,n)$,
\begin{equation}
\begin{array}{cccccccccccccccccc}
\pi_2 & \pi_4 &  \cdots & \pi_{2i} &\pi_{2i+2} & \cdots& \pi_{2n}\\[8pt]
\pi_1 & \pi_3 &  \cdots & \pi_{2i-1} &\pi_{2i+1} & \cdots & \pi_{2n-1}.
\end{array}
\end{equation}
As a result, by  formula \eqref{formula}, we obtain
\begin{align*}
f_n(2n)=f^{(n,n)}&=(2n)!
\begin{vmatrix}
\dfrac{1}{n!} & \dfrac{1}{(n+1)!} \\[12pt]
\dfrac{1}{(n-1)!}   & \dfrac{1}{n!}
\end{vmatrix}\\[12pt]
&=\frac{1}{n+1}{2n\choose n}\\[12pt]
&=C_n.
\end{align*}

So far we have established a one-to-one correspondence between the
set $\mathscr{F}_{a_1,a_2,\ldots,a_k}$ and the set of
$2$-regular skew tableaux in Theorem \ref{f}. Hence the
enumeration of the number of minimal permutations is
equivalent to the enumeration of skew Young tableaux.
\begin{corollary}\label{g}
Given an ascent sequence $\alpha=(a_1,a_2,\ldots,a_k)$, where
$\sum_{i=1}^ka_i=n$ and $a_i\geq 2$, for $1\leq i\leq k$, we have
\begin{equation}\label{deter}
F_{a_1,a_2,\ldots,a_k}=n!\det(A)=n! \begin{vmatrix}
\dfrac{1}{a_1!}\\[8pt]
1           &\dfrac{1}{a_2!}& \multicolumn{3}{c}{\raisebox{1.3ex}[0pt]{$A_{ij}$}}\\[8pt]
  A_{3,1}       & 1            &\dfrac{1}{a_3!}\\[8pt]
  &   A_{4,2}  & 1   &\dfrac{1}{a_4!}\\[8pt]
&  &  \ddots & \ddots &\ddots \\[8pt]
&&& A_{k-1,k-3} & 1 & \dfrac{1}{a_{k-1}!}\\[8pt]
&&&&A_{k,k-2}&1&\dfrac{1}{a_k!}
\end{vmatrix},
\end{equation}
where
\begin{align*}
A_{i,j}&=0, \quad\quad\text{when}\quad j<i-2,\\[6pt]
A_{i,i-2}&=\begin{cases}
0,\quad\quad\quad\quad \mbox {if}\quad a_{i-1}> 2,\\[3pt]
1,\quad\quad\quad\quad \mbox {if}\quad a_{i-1}=2,
\end{cases}\\[6pt]
A_{i,j}&=\dfrac{1}{\left(\sum_{m=i}^ja_m-(j-i)\right)!},\quad\quad \text{when}\quad j>i.
\end{align*}
\end{corollary}

\pf
First of all we need to determine the shape of the $2$-regular skew tableau $P$ defined in
\eqref{aa}.
Suppose the shape of $P$ is $\lambda/\mu$, by the
correspondence described in Theorem \ref{f}, the number of elements in
each column of $P$ are $a_1, a_2,\ldots$ and $a_k$,
respectively.
In other words,
\begin{equation}\label{1}
\lambda_i'-\mu_i'=a_i,\quad \mbox{for}\quad 1\leq i\leq k.
\end{equation}
Furthermore, the fact that the uppermost two elements in the
$i$th column and the lowest two elements in the $(i+1)$th
column compose a $2\times 2$ square leads to
\begin{equation}\label{2}
\lambda_{i}'-\lambda_{i+1}'=a_{i}-2, \quad \mbox{for} \quad 1\leq i\leq k-1.
\end{equation}
Obviously $\lambda'_k=a_k$ and $\mu'_k=0$. For
$1\leq i\leq k-1$, we have
\begin{align*}
\lambda'_{i}&=\lambda_{i+1}'+(a_{i}-2)\\[6pt]
&=\lambda'_{i+2}+(a_{i+1}-2)+(a_{i}-2)\\[6pt]
&=\cdots\\[6pt]
&=\lambda'_k+(a_{k-1}-2)+\cdots+(a_{i}-2)\\[6pt]
&=a_k+(a_{k-1}-2)+\cdots+(a_{i}-2)\\[6pt]
&=a_{i}+a_{i+1}+\cdots+a_k-2(k-i).
\end{align*}
This yields that $\mu'_{i}=\lambda'_{i}-a_{i}
=a_{i+1}+\cdots+a_k-2(k-i)$.
So $\lambda/\mu$ is exactly the shape of $P$.
As a consequence,
\begin{equation}\label{mai}
F_{a_1,a_2,\ldots,a_k}(n)=f^{\lambda'/\mu'},
\end{equation}
where
$$\lambda'=(\lambda'_1,\lambda'_2,\cdots,\lambda'_k),
\quad \lambda'_i=\sum_{j=i}^ka_j-2(k-i),
\quad \mbox{for}\quad 1\leq i\leq k,$$
and
$$\mu'=(\mu'_1,\mu'_2\ldots,\mu'_k),
\quad \mu'_i=\lambda'_i-a_i=\sum_{j=i+1}^ka_j-2(k-i), \quad \mbox{for}\quad 1\leq i\leq k.$$
Now we proceed to compute the number of $2$-regular skew tableaux
by using formula \eqref{formula}.
We consider the shape $\lambda'/\mu'$
of the $2$-regular skew tableau $P$ by
dividing the relations between $i$ and $j$
into the following  cases.
From the equations \eqref{1} and \eqref{2}, we obtain
\begin{itemize}
\item[(1).]If $j< i-2$,
$\lambda'_i-\mu_{j}'=-(a_{j+1}+\cdots+a_{i-1})+2(i-j)$.
Therefore, $\lambda'_i-\mu_{j}'-i+j
=-(a_{j+1}+\cdots+a_{i-1})+(i-j)$.
Since  $a_m\geq 2 (1\leq m\leq k)$, we see that
 $\lambda'_i-\mu_{j}'-i+j<0$,
 which means $A_{i,j}=0$.
\item[(2).]If $j=i-2$, $\lambda'_i-\mu_{i-2}'
=-a_{i-1}+4$. Therefore, if $a_{i-1}=2$,
then $\lambda'_i-\mu_{i-2}'-i+{i-2}=0$
and $A_{i,i-2}=1/0!=1$. Otherwise, we have
 $a_{i-1}>2$, $\lambda'_i-\mu_{i-2}'-i+{i-2}<0$,
 and $A_{i,i-2}=0$.
\item[(3).]If $j=i-1$, $\lambda'_i-\mu_{i-1}'=2$.
Consequently, $\lambda'_i-\mu_{i-1}'-i+(i-1)=1$,
and $A_{i,i-1}=1$.
\item[(4).]If $j=i$, then from the relation \eqref{1}, we have
 $\lambda'_i-\mu_i'=a_i$, and so $A_{i,i}=\dfrac{1}{a_i!}$.
\item[(5).]If $j\geq i-1$,
then by \eqref{2}, $\lambda'_i-\mu'_j
=a_i+(a_{i+1}-2)+\cdots+(a_j-2)
=a_i+\cdots a_j-2(j-i)$. In this case,
    $$A_{i,j}=\frac{1}{(\lambda'_i-\mu'_j-i+j)!}
    =\frac{1}{\left(a_i+\cdots+a_j-(j-i)\right)!},$$
\end{itemize}
as desired. This completes the proof.
\qed

We remark that  Corollary \ref{g}
 can be viewed as a refinement of the number $f_d(n)$,
 that is the number of minimal permutations
 in $\mathcal{F}_{n-k}(n)$ with prescribed ascent set
 $\{a_1,a_1+a_2,\ldots,a_1+\cdots+a_{k-1}\}$.

We can now compute the number of minimal
permutations in $\mathcal{F}_d(n)$. Note that such
 minimal permutations have $n-d$ maximal
decreasing substrings.
\begin{corollary}\label{main2}
For $d+1\leq n\leq 2d$, we have
\begin{equation}\label{fa}
f_{d}(n)=\sum_{\begin{subarray}{l}a_i\geq 2\,
\text{for}\,1\leq i\leq n-d\\a_1+a_2+\cdots+a_{n-d}
=n\end{subarray}}F_{a_1,a_2,\ldots,a_{n-d}}.
\end{equation}
\end{corollary}

As an application of the above
formula \eqref{fa} for $d=n-2$,
we immediately obtain the formula for $f_{n-2}(n)$
due to Bouvel and Pergola \cite{bou}.
It is obvious that the minimal permutations
 in $\mathcal{F}_{n-2}(n)$ have only one ascent,
 which implies that they have two maximal decreasing
  substrings. Suppose that the unique ascent is $k$,
  and the ascent sequence is $(k,n-k)$
  for $2\leq k\leq n-2$. By Corollary \ref{g},
  it is easy to check
\[
F_{k,n-k}=n!
\begin{vmatrix}
\dfrac{1}{k!} & \dfrac{1}{(n-1)!}\\[12pt]
     1       & \dfrac{1}{(n-k)!}
\end{vmatrix}={n\choose k}-n.
\]
By Corollary \ref{main2}, we arrive at
$$f_{n-2}(n)=
\sum_{k=2}^{n-2}\left({n\choose k}-n\right)=2^n-2-n(n-1).$$
We now come to the computation of  $f_{n-3}(n)$.

\begin{theorem}
The number of minimal permutations of
length $n$ with $n-3$ descents equals
$$f_{n-3}(n)=3^n-(n^2-2n+4)
2^{n-1}+\frac{1}{2}\left(n^4-7n^3+19n^2-21n+2\right).$$
\end{theorem}

\pf The minimal permutations in $\mathcal{F}_{n-3}(n)$
have three maximal decreasing substrings. Suppose the ascent
 sequence is $(a,b,c)$ such that $a+b+c=n$ and $a,b,c\geq 2$.
 According to Corollary \ref{g}, let
\begin{align*}
A_1&=n!
\begin{vmatrix}
\dfrac{1}{a!} & \dfrac{1}{(a+1)!}  & \dfrac{1}{(n-2)!}\\[12pt]
    1        & \dfrac{1}{2!}     &  \dfrac{1}{(c+1)!} \\[12pt]
    1        &     1            &  \dfrac{1}{c!}
\end{vmatrix}\\[8pt]
&=\frac{n!}{a!2!c!}+\frac{n!}{(a+1)!(c+1)!}+\frac{n!}{(n-2)!}\\[8pt]
& \qquad-\frac{n!}{(n-2)!2!}
-\frac{n!}{(a+1)!c!}-\frac{n!}{a!(c+1)!},
\end{align*}
and let
\begin{align*}
A_2&=n!
\begin{vmatrix}
\dfrac{1}{a!} & \dfrac{1}{(a+b-1)!}  & \dfrac{1}{(n-2)!}\\[12pt]
    1        & \dfrac{1}{b!}     &  \dfrac{1}{(b+c-1)!} \\[12pt]
    0        &     1            &  \dfrac{1}{c!}
\end{vmatrix}\\[8pt]
&=\frac{n!}{a!b!c!}+\frac{n!}{(n-2)!}
-\frac{n!}{(a+b-1)!c!}-\frac{n!}{a!(b+c-1)!}.
\end{align*}
By Corollary \ref{main2}, we obtain
\[
f_{n-3}(n)
=\sum_{\begin{subarray}{c}a,c\geq 2,b=2\\a+b+c=n\end{subarray}}A_1+\sum_{\begin{subarray}{c}a,c\geq 2,b\geq3\\a+b+c=n\end{subarray}}A_2.
\]
In order to simplify the computation,
we reformulate the above equation into the following form,
\begin{align}\label{n-3}
f_{n-3}(n)&=\sum_{\begin{subarray}{c}a,c\geq 2,b=2\\a+b+c=n\end{subarray}}A_1+\sum_{\begin{subarray}{c}a,b,c\geq 2\\a+b+c=n\end{subarray}}A_2-\sum_{\begin{subarray}{c}a,c\geq 2,b=2\\a+b+c=n\end{subarray}}A_2 \nonumber \\[6pt]
&=\sum_{\begin{subarray}{c}a,c\geq 2,b=2\\a+b+c=n\end{subarray}}(A_1-A_2)+\sum_{\begin{subarray}{c}a,b,c\geq 2\\a+b+c=n\end{subarray}}A_2.
\end{align}
It is easy to check that
\begin{align}\label{w}
\sum_{\begin{subarray}{c}a,c\geq 2,b=2\\a+b+c=n\end{subarray}}(A_1-A_2)&
=\sum_{\begin{subarray}{c}a,c\geq 2\\a+c=n-2\end{subarray}}\left(\frac{n!}{(a+1)!(c+1)!}-{n\choose 2}\right)\nonumber\\[6pt]
&=\sum_{a=2}^{n-4}\left({n\choose a+1}-{n\choose 2}\right) \nonumber\\[6pt]
&=2^n-2n-2-(n-3){n\choose 2}.
\end{align}
The second sum of \eqref{n-3} can be expressed as follows
$$\sum_{\begin{subarray}{c}a,b,c\geq 2\\a+b+c=n\end{subarray}}A_2=\sum_{\begin{subarray}{c}a,b,c\geq 2\\a+b+c=n\end{subarray}}\left({n\choose a,b,c}+n(n-1)-n{n-1\choose c}-n{n-1\choose a}\right).$$
On the one hand, by the inclusion-exclusion principle, we have
\begin{align}\label{v1}
\sum_{\begin{subarray}{c}a,b,c\geq 2\\a+b+c=n\end{subarray}}{n\choose a,b,c}&=3^n-3\sum_{b=0}^n{n\choose 0,b,n-b}-3\sum_{b=0}^{n-1}{n\choose 1,b,n-1-b} \nonumber \\[6pt]
&\qquad+3{n\choose 0,0,n}+3{n\choose 1,1,n-2}+6{n\choose 0,1,n-1}\nonumber \\[6pt]
&=3^n-3\cdot 2^n-3n\cdot 2^{n-1}+3n^2+3n+3.
\end{align}
On the other hand, let $[x^n]f(x)$ denote the coefficient of
$x^n$ in $f(x)$, then we get
\begin{align}\label{v2}
\sum_{\begin{subarray}{c}a,b,c\geq 2\\a+b+c=n\end{subarray}}
n(n-1)
 &=n(n-1)\cdot[x^n](x^2+x^3+\cdots)^3 \nonumber \\
&=n(n-1){n-4\choose 2}.
\end{align}
Furthermore,
\[
 n\sum_{\begin{subarray}{c}a,b,c\geq 2\\a+b+c=n\end{subarray}}
 \left({n-1\choose c}+{n-1\choose a}\right)=
2n\sum_{\begin{subarray}{c}a,b,c\geq 2\\a+b+c=n\end{subarray}}{n-1\choose a}\]
It is easily seen that
\[
\sum_{\begin{subarray}{c}a,b,c\geq 2\\a+b+c=n\end{subarray}}
{n-1\choose a}=
\sum_{a=2}^{n-4} \sum_{\begin{subarray}{c}b,c\geq 2\\b+c=n-a\end{subarray}} {n-1\choose a}  =
\sum_{a=2}^{n-4}{n-1\choose a}(n-3-a).
\]
Since $$\sum_{a=1}^{n-1}a{n-1\choose a}
=\sum_{a=1}^{n-1}a{n-1\choose a}x^{a-1}\Bigg |_{x=1}
=(n-1)(1+x)^{n-2}\bigg|_{x=1}=(n-1)2^{n-2},$$
we find
\begin{align}\label{v3}
2n\sum_{\begin{subarray}{c}a,b,c\geq 2\\a+b+c=n\end{subarray}}{n-1\choose a}&=2n(n-3)\sum_{a=2}^{n-4}{n-1\choose a}-2n\sum_{a=2}^{n-4}a{n-1\choose a} \nonumber \\[6pt]
&=2n(n-3)\left(2^{n-1}-2(n-1)-2\right) \nonumber \\[6pt]
&\qquad -2n\left((n-1)2^{n-2}-2(n-1)-(n-1)(n-2)\right) \nonumber\\[6pt]
&=n(n-3)2^{n}-n(n-1)2^{n-1}-2n^3+10n^2.
\end{align}
By \eqref{w}, \eqref{v1}, \eqref{v2} and \eqref{v3}, we finally
obtain
\[
f_{n-3}(n)
=3^n-(n^2-2n+4)2^{n-1}+\frac{1}{2}\left(n^4-7n^3+19n^2-21n+2\right),
\]
as claimed.
\qed

\section{A refinement of $f_{n+1}(2n+1)$ via Knuth equivalence}

In this section, we give a combinatorial proof of a
refined formula for the number $f_{n+1}(2n+1)$.
Given a minimal permutation $\pi=
\pi_1\pi_2\cdots\pi_{2n+1}$ of length $2n+1$ with $n+1$ descents,
 there are only one occurrence of consecutive descents in $\pi$ and the
 other descents are separated by ascents.
  We shall consider the set of
  minimal permutations for which the
  unique consecutive descents are $2i-1$
  and $2i$.

\begin{theorem}\label{mi}
Let $\mathcal{M}_{2n+1,2i}$ be the subset of $\mathcal{F}_{n+1}(2n+1)$ whose unique consecutive descents are $2i-1$ and $2i$, where $1\leq i\leq n$. We have
\begin{equation}\label{refine}
|\mathcal{M}_{2n+1,2i}|={2n+1\choose n-1}{n-1\choose i-1}.
\end{equation}
\end{theorem}
We shall  give two proofs of this theorem.

It is easy to see that the $2$-regular skew
tableaux corresponding to minimal permutations
$\pi\in\mathcal{M}_{2n+1,2i}$  are of the following form,
\begin{equation}\label{skew}
\begin{array}{ccccccccccccccccc}
&&&&&& \pi_{2i+1} & & \pi_{2i+3} &  & \pi_{2i+5} &  & \cdots &  & \pi_{2n-1}& & \pi_{2n+1}\\
\pi_2 &  & \pi_4 &  & \cdots &  & \pi_{2i} &  &\pi_{2i+2} & & \pi_{2i+4} & & \cdots & & \pi_{2n-2} & & \pi_{2n}\\
\pi_1 & & \pi_3 && \cdots && \pi_{2i-1}.
\end{array}
\end{equation}
The conjugate shape $\lambda'/\mu'$ is
$$(\underbrace{3,3,\ldots,3}_i,\underbrace{2,2,\ldots,2}_{n-i})
/(\underbrace{1,1,\ldots,1}_{i-1}),
\quad \mbox{for} \quad 1\leq i\leq n.$$
Notice that the skew shape $\lambda/\mu$ is $(n,n,i)/(i-1)$.
In this context, we can obtain the number of
 skew Young tableaux directly from formula \ref{formula},
\begin{align*}
f^{(n,n,i)/(i-1)}&=(2n+1)!
\begin{vmatrix}
\dfrac{1}{(n-i+1)!} & \dfrac{1}{(n+1)!}  &  \dfrac{1}{(n+2)!} \\[12pt]
\dfrac{1}{(n-i)!}   & \dfrac{1}{n!}      & \dfrac{1}{(n+1)!}\\[12pt]
  0                & \dfrac{1}{(i-1)!}  & \dfrac{1}{i!}
\end{vmatrix} \\[12pt]
&={2n+1\choose n}{n+1\choose i}+{2n+1\choose n-1}{n-1\choose i-1}\\[12pt]
&\qquad-{2n+1\choose n}{n\choose i}-{2n+1\choose n}{n\choose i-1}\\[12pt]
&={2n+1\choose n-1}{n-1\choose i-1}.
\end{align*}
So we immediately get
$$f_{n+1}(2n+1)=\sum_{i=1}^n f^{(n,n,i)/(i-1)}=2^{n-1}{2n+1\choose n-1}.
$$
\qed

We next give a combinatorial interpretation of Theorem \ref{mi}.
We give an overview of the background on
the RSK correspondence and  the Knuth equivalence.

 Suppose $\pi\stackrel{\mathrm{RSK}}{\longrightarrow}(P,Q)$,
 where $P$ is called the insertion tableau while $Q$
 is called the recording tableau. Two
permutations are Knuth-equivalent if and only if
their insertion tableaux are the same.

The following properties of insertion paths will be useful.
Denote by $I(P\leftarrow k)$ the
insertion path of a positive integer $k$
into an SYT (standard Young tableau)
 $P=(P_{ij})$ by the RSK algorithm. Then
\begin{enumerate}
\item[(a)] When we insert $k$ into an SYT $P$, the insertion path moves to the left. More precisely, if $(r,s),(r+1,t)\in I(P\leftarrow k)$ then $t\leq s$.
\item[(b)] Let $P$ be an SYT, and let $j< k$. Then $I(P\leftarrow j)$ lies strictly to the left of $I((P\leftarrow j)\leftarrow k)$. More precisely, if $(r,s)\in I(P\leftarrow j)$, and $(r,t)\in I((P\leftarrow j)\leftarrow k)$, then $s<t$. Moreover, $I((P\leftarrow j)\leftarrow k)$ does not extend below the bottom of $I(P\leftarrow j)$. Equivalently,
    $$\# I((P\leftarrow j)\leftarrow k) \leq\# I(P\leftarrow j).$$
\end{enumerate}
See \cite{sta} for more details.

Now we begin the combinatorial proof by using the Knuth equivalence of permutations in connection
with the RSK correspondence
between permutations and standard Young tableaux, see Stanley
\cite{sta}.

Let $\mathcal{T}_{2n+1,k}$ be
 the set of standard Young tableaux of
  size $2n+1$ with shape $(n,n+1-k,k)$,
   where $1\leq k\leq [\frac{n+1}{2}]$,
   where $[x]$ denotes the largest integer not exceeding  $x$.
 According to the correspondence between
 minimal permutations and  $2$-regular skew tableaux,
  we see that $\mathcal{M}_{2n+1,2i}$ is
   the set of the $2$-regular skew tableaux
   with row lengths $n-i+1,n,i$.

\begin{theorem}\label{main}
There exists a bijection between $\mathcal{M}_{2n+1,2i}$
 and $\mathcal{T}_{2n+1,k}$ such that the
length of the last row of an SYT
in $\mathcal{T}_{2n+1,k}$ does not exceed
the smallest row length of the corresponding
$2$-regular tableau in $\mathcal{M}_{2n+1,2i}$.
Equivalently, we have the following formula,
\begin{equation}\label{symm}
|\mathcal{M}_{2n+1,2i}|=\sum_{k=1}^j|\mathcal{T}_{2n+1,k}|,\quad\quad \mbox{where}\quad j=\min\left\{n-i+1,i\right\}.
\end{equation}
Consequently, $|\mathcal{M}_{2n+1,2i}|$ is symmetric in $i$,
\begin{equation}\label{m}
|\mathcal{M}_{2n+1,2i}|=|\mathcal{M}_{2n+1,2(n-i+1)}|.
\end{equation}
\end{theorem}

The main idea of the proof can be described as follows.
For every $\pi\in\mathcal{M}_{2n+1,2i}$,
we show that there always exists a permutation
$\pi'\in\mathfrak{S}_{2n+1}$ which equivalent to $\pi$
(by the Knuth equivalence, to be precise). In other words,
they have the same insertion tableau.
Therefore, by constructing the insertion tableau of $\pi'$,
we can give a description of the insertion tableau of $\pi$.
Thus we obtain the shapes of the SYTs
corresponding to the permutations in $\mathcal{M}_{2n+1,2i}$.

\noindent
{\it Proof.}
Let $
\pi'=\pi_1\pi_2\cdots\pi_{2i-1}\pi_{2i}\pi_{2i+2}
 \cdots\pi_{2n}\pi_{2i+1}\pi_{2i+3}\cdots\pi_{2n+1}\in\mathfrak{S}_{2n+1},$
that is, $\pi'$ is obtained from
$\pi$ by fixing the first $2i$ elements,
and moving the remaining elements with even subscripts
forward and those with odd subscripts backward.
By \eqref{skew}, $\pi'$ can also obtained by first
reading the elements of the last two rows of
\eqref{skew} and keeping  the order of these
elements in $\pi$ unchanged.
Then read off the elements of the first row of \eqref{skew}.

First, we show that $\pi$ and $\pi'$ are Knuth equivalent, namely,
$$\pi\stackrel{K}{\thicksim}\pi'.$$
Recall that each Knuth
transformation switches two adjacent entries $a$ and $c$
provided that an entry $b$ satisfying $a<b<c$ is
located next to $a$ or $c$.
Write
$$\pi=\pi_1\pi_2\cdots\pi_{2i-1}
\pi_{2i}\boxed{\pi_{2i+1}}\pi_{2i+2}
\boxed{\pi_{2i+3}}\cdots\pi_{2n}\boxed{\pi_{2n+1}}.$$
For the purpose of presentation,
 the elements which will be moved back are
  framed. When we write $bac$ under three  consecutive elements,
   we mean that these three elements have type $bac$.
  We shall apply a series of Knuth transformations
  to the substring $\pi_{2i}\boxed{\pi_{2i+1}}
  \pi_{2i+2}\cdots\boxed{\pi_{2n+1}}$ of $\pi$.
  This is equivalent
to exchanging every two adjacent elements after $\pi_{2i}$,
\begin{align*}
\pi &\, =\pi_1\cdots \pi_{2i-1}\underbrace{\pi_{2i}\boxed{\pi_{2i+1}}\pi_{2i+2}}_{bac}
          \underbrace{\boxed{\pi_{2i+3}}\pi_{2i+4}\boxed{\pi_{2i+5}}}_{acb}\pi_{2i+6}\cdots \boxed{\pi_{2n-1}}
           \pi_{2n}\boxed{\pi_{2n+1}}\\[8pt]
&\stackrel{\rm{K}}\thicksim\pi_1 \cdots \pi_{2i-1}\pi_{2i}\pi_{2i+2}\boxed{\pi_{2i+1}}\pi_{2i+4}\boxed{\pi_{2i+3}}
         \underbrace{\boxed{\pi_{2i+5}}\pi_{2i+6}\boxed{\pi_{2i+7}}}_{acb} \pi_{2i+8}\boxed{\pi_{2i+9}}\pi_{2i+10} \boxed{\pi_{2i+11}} \cdots  \\
&\quad\cdots\cdots\\
&\stackrel{\rm{K}}\thicksim \pi_1\cdots\pi_{2i-1} \pi_{2i} \underbrace{\pi_{2i+2} \boxed{\pi_{2i+1}} \pi_{2i+4}}_{bac}\underbrace{\boxed{\pi_{2i+3}}\pi_{2i+6}
         \boxed{\pi_{2i+5}}}_{acb} \pi_{2i+8}\boxed{\pi_{2i+9}} \cdots
        \pi_{2n}\boxed{\pi_{2n-1}}\boxed{\pi_{2n+1}}.
\end{align*}
By this procedure we have moved $\pi_{2i+2}$
forward and $\pi_{2n-1}$ backward,
Then for the substring
$\pi_{2i+2}\boxed{\pi_{2i+1}}\pi_{2i+4}\cdots\boxed{\pi_{2n+1}}$
of the resulting permutation, repeat the above procedure in order to
move $\pi_{2i+4}$ forward and $\pi_{2n-3}$ backward.
\begin{align*}
\pi &\, \stackrel{\rm{K}}\thicksim \pi_1\cdots\pi_{2i}\pi_{2i+2}\pi_{2i+4} \boxed{\pi_{2i+1}} \pi_{2i+6}\boxed{\pi_{2i+3}}
         \underbrace{\boxed{\pi_{2i+5}}\pi_{2i+8}\boxed{\pi_{2i+7}}}_{acb}\pi_{2i+10} \boxed{\pi_{2i+9}}\cdots\\[8pt]
&\stackrel{\rm{K}}\thicksim\pi_1\cdots\pi_{2i}\pi_{2i+2}\pi_{2i+4}\boxed{\pi_{2i+1}}\pi_{2i+6}\boxed{\pi_{2i+3}}\pi_{2i+8}
         \boxed{\pi_{2i+5}}\underbrace{\boxed{\pi_{2i+7}}\pi_{2i+10}\boxed{\pi_{2i+9}}}_{acb}\cdots\\
         &\quad\cdots\cdots\\
&\stackrel{\rm{K}}\thicksim\pi_1\cdots\pi_{2i}\pi_{2i+2}\underbrace{\pi_{2i+4}\boxed{\pi_{2i+1}}\pi_{2i+6}}_{bac}\underbrace{\boxed{\pi_{2i+3}}
         \pi_{2i+8}\boxed{\pi_{2i+5}}}_{acb}\pi_{2i+10}\cdots
         \pi_{2n}\boxed{\pi_{2n-3}}\boxed{\pi_{2n-1}}\boxed{\pi_{2n+1}}.
\end{align*}
Iterating this process until all the elements
after $\pi_{2i}$ with even subscripts are
moved forward while the elements
after $\pi_{2i}$ of odd subscripts are moved backward, we get
\begin{align*}
\pi &\, \stackrel{\rm{K}}\thicksim\pi_1\cdots\pi_{2i-1}\pi_{2i}\pi_{2i+2}\pi_{2i+4}\cdots\pi_{2n-2}\pi_{2n}\boxed{\pi_{2i+1}}\boxed{\pi_{2i+3}}
\boxed{\pi_{2i+5}}\cdots\boxed{\pi_{2n-3}}\boxed{\pi_{2n-1}}\boxed{\pi_{2n+1}}\\[12pt]
&=\pi'.
\end{align*}
We give an example  to illustrate the above procedure.
Let $\pi=6\,3\,7\,4\,1\,5\,2\,9\,8\,11\,
10\,13\,12\in \mathcal{F}_7(13)$. We have the following Knuth
transformations:
\begin{align*}
\pi&= 6\,3\,7\,4\,\boxed{1}\,5\,\boxed{2}\,9\,\boxed{8}\,11\,\boxed{10}\,13\,\boxed{12}
\stackrel{\rm{K}}\thicksim 6\,3\,7\,4\,5\,\boxed{1}\,9\,\boxed{2}\,\boxed{8}\,11\,\boxed{10}\,13\,\boxed{12}\\[12pt]
&\ \stackrel{\rm{K}}\thicksim 6\,3\,7\,4\,5\,\boxed{1}\,9\,\boxed{2}\,11\,\boxed{8}\,\boxed{10}\,13\,\boxed{12}
 \stackrel{\rm{K}}\thicksim 6\,3\,7\,4\,5\,\boxed{1}\,9\,\boxed{2}\,11\,\boxed{8}\,13\,\boxed{10}\,\boxed{12}\\[12pt]
&\ \stackrel{\rm{K}}\thicksim 6\,3\,7\,4\,5\,9\,\boxed{1}\,11\,\boxed{2}\,\boxed{8}\,13\,\boxed{10}\,\boxed{12}
\stackrel{\rm{K}}\thicksim 6\,3\,7\,4\,5\,9\,\boxed{1}\,11\,\boxed{2}\,13\,\boxed{8}\,\boxed{10}\,\boxed{12}\\[12pt]
&\ \stackrel{\rm{K}}\thicksim 6\,3\,7\,4\,5\,9\,11\,\boxed{1}\,\boxed{2}\,13\,\boxed{8}\,\boxed{10}\,\boxed{12}
 \stackrel{\rm{K}}\thicksim 6\,3\,7\,4\,5\,9\,11\,\boxed{1}\,13\,\boxed{2}\,\boxed{8}\,\boxed{10}\,\boxed{12}\\[12pt]
&\ \stackrel{\rm{K}}\thicksim 6\,3\,7\,4\,5\,9\,11\,13\,\boxed{1}\,\boxed{2}\,\boxed{8}\,\boxed{10}\,\boxed{12}\\[12pt]
&=\pi'.
\end{align*}
Therefore,
$\pi$ and $\pi'$ have the same insertion tableau.
Applying the RSK algorithm to $\pi'$, it is easy to see that the standard Young tableau corresponding to the first $n+i$ elements $\pi_1\pi_2\cdots\pi_{2i}\pi_{2i+2}\cdots\pi_{2n}$ of $\pi'$ is an SYT of shape $(n,i)$,
\begin{equation}\label{p'}
P'=\begin{array}{ccccccc}
\pi_2 & \pi_4 & \cdots & \pi_{2i} & \pi_{2i+2} & \cdots & \pi_{2n}\\
\pi_1 & \pi_3 & \cdots & \pi_{2i-1},
\end{array}
\end{equation}
which is exactly the last two rows of \eqref{skew}.
The insertion tableau of $\pi'$ can be obtained as follows
\begin{equation}\label{insert}
\bigg(\big(\cdots\left(\left(P'\leftarrow\pi_{2i+1}\right)\leftarrow\pi_{2i+3}\right)\cdots\leftarrow \big)\pi_{2n-1}\bigg)\leftarrow\pi_{2n+1}.
\end{equation}
Since $\pi$ and $\pi'$ have
 the same insertion tableau, \eqref{insert}
can be also considered as the insertion tableau of $\pi$.

We now aim to give a second combinatorial proof.
First we show that
$$\mathcal{M}_{2n+1,2i}\longrightarrow \bigcup_{k=1}^j
 \mathcal{T}_{2n+1,k}$$ is an injection.
Since $\pi_{2i+1}>\pi_{2i}>\pi_{2i-1}$,
when inserting $\pi_{2i+1}$ into $P'=(P_{ij})$,
 the resulting tableau corresponding $P_0=P'\leftarrow \pi_{2i+1}$
  is of shape $(n,i,1)$. Moreover, the intersection position
   of $I(P_0)$ in the first row cannot be to the
    right of $(i,1)$ in $P'$.
By induction, we assume that the
insertion tableau of $P_{m-1}=
(\cdots(P'\leftarrow \pi_{2i+1})\cdots\leftarrow
\pi_{2i+2(m-1)+1})$ is of shape $(n,i+m-s,s)$,
where $1\leq s\leq m$. Then let us examine the
shape of $P_m=P_{m-1}\leftarrow \pi_{2i+2m+1}$.
Since $\pi_{2i+2m+1}>\pi_{2i+2m-1}$, the insertion
path $I(P_m)$ lies strictly to the right of
$I(P_{m-1})$ and does not extend below the bottom
 of $I(P_{m-1})$.
 Since $\pi_{2i+2m+1}>\pi_{2i+2m}$, the insertion
 path of $P_m$ in the first row cannot extend to
 the right of $(i+m,1)$ of $P'$.
 It follows that the shape of $P_m$ can be
  obtained from that of  $P_{m-1}$ by adding a
  new element to the second or the third row
  of $P_{m-1}$. We deduce
   that the shape of $P_m$ must be one of the form
   $(n,i+m+1-s,s)$, where $1\leq s\leq m+1$.

 Next we aim to show  that
  $$\bigcup_{k=1}^j
  \mathcal{T}_{2n+1,k}
  \longrightarrow\mathcal{M}_{2n+1,2i}$$
  is also an injection.
  Given a standard Young tableau
  $P$ of shape $(n,n+1-k,k)$.
  Pick up the set of positions
  in $P$ which are note occupied
   by elements in $P'$ given
   by \eqref{p'}. Suppose that
   these positions are $(i_1,j_1),(i_2,j_2),\ldots,
   (i_{n-i},j_{n-i})$ such that
    $j_1\geq j_2\cdots \geq j_{n-i}$
    and that  if $j_t=j_{t+1}$,
    then $i_t<i_{t+1}$.
    In other words, these positions are
    ordered from the northeast corner to the southwest corner.

At first,
we apply  the ``inverse bumping''
to $P_{i_1j_1}$. It bumps an element
$\pi_{r_1t_1}$ in the first row
from $P$. Put $\pi_{r_1t_1}$ on top of $\pi_{1n}$
of $P\rightarrow P_{i_1j_1}$.
It is easy to see that $\pi_{r_1t_1}<\pi_{1n}$.
Note that when we begin to apply the ``inverse bumping"
to $P_{i_1j_1}$, it is put the end of its
row (row $i_1$). Inductively,
suppose $P_{i_sj_s}$ bumps some element
 $\pi_{r_st_s}$ in the first row of $P$,
  and $\pi_{r_st_s}$ is putted on top
  of $\pi_{1,n-s+1}$, where $\pi_{r_st_s}<\pi_{1,n-s+1}$.
   When we apply the inverse bumping
   to $P_{i_{s+1}j_{s+1}}$, its ``inverse insertion path"
   intersecting row $i_s$ is strictly to
   the left of column $j_s$.
   Consequently, at  row $i_s$,
    the inverse insertion path
    of $P_{i_{s+1}j_{s+1}}$ lies
    strictly to the left of that of $P_{i_sj_s}$.
    By induction, the entire inverse insertion path
     of $P_{i_{s+1}j_{s+1}}$ lies strictly to
      the left of that of $P_{i_sj_s}$.
       In particular, the element $\pi_{r_{s+1}t_{s+1}}$
        bumped by $P_{i_{s+1}j_{s+1}}$ in the first row
         is to the left of $\pi_{r_st_s}$.
         Hence  $\pi_{r_{s+1}t_{s+1}}<\pi_{r_st_s}$.

We now put $\pi_{r_{s+1}t_{s+1}}$ to the left
 of $\pi_{r_st_s}$. Since the inverse insertion path
  of $P_{i_{s+1}j_{s+1}}$ lies strictly to the left of
   that of $P_{i_sj_s}$, we find $\pi_{r_{s+1}t_{s+1}}<\pi_{1\,n-s}$,
   as required.
Note that the condition
$j=\min{i,n+1-i}$ is necessary
since  the
resulting tableau is a standard Young tableau.

The symmetry of $\mathcal{M}_{2n+1,2i}$ is immediate from \eqref{symm}.
\qed

We now aim to compute the number of SYTs
of shape $(n,n+1-k,k)$.
Recall that
if $\lambda\vdash n$, then the number of SYTs of shape $\lambda$ is given by the hook length formula,
$$f^{\lambda}=\frac{n!}{\prod_{u\in\lambda}h(u)},$$
where $h(u)=\lambda_i+\lambda_j'-i-j+1.$

By the hook length formula, it is easy to show that ${2n+1\choose n-1}$ counts the number of SYTs of shape $(n,n,1)$. For the general case, we have
\begin{theorem}\label{syt}
For $2\leq k\leq [\frac{n+1}{2}]$, the number of SYTs of shape $(n,n+1-k,k)$ is
\begin{equation}\label{hook}
|\mathcal{T}_{2n+1,k}|=\frac{n-2k+2}{k-1}{n-1\choose k-2}{2n+1\choose n-1}.
\end{equation}
\end{theorem}

\proof We conduct induction on $k$.
When $k=2$, by the hook length formula,
 it is easy to show that the number of SYTs of
 shape $(n,n-1,2)$ equals $$(n-2){2n+1\choose n-1}.$$
 We now suppose that \eqref{hook}
 holds for $k-1$. Comparing the hook lengths of SYTs of
 shape $(n,n-k+1,k)$ to those of shape $(n,n+2-k,k-1)$,
 we find that
$$\frac{|\mathcal{T}_{2n+1,k-1}|}{|\mathcal{T}_{2n+1,k}|}
=\frac{(n-2k+2)(n-k+2)}{(n-2k+4)(k-1)}.$$
Thus the number of SYTs of shape $(n,n-k+1,k)$ is given by
\begin{align*}
|\mathcal{T}_{2n+1,k}|&=\frac{(n-2k+2)(n-k+2)}{(n-2k+4)(k-1)}\frac{(n-2k+4)}{k-2}{n-1\choose k-3}{2n+1\choose n-1}\\[8pt]
&=\frac{n-2k+2}{k-1}{n-1\choose k-2}{2n+1\choose n-1}.
\end{align*}
 \qed

We are now ready to complete the proof of Theorem \ref{mi}.

\noindent{\it The proof of Theorem \ref{mi}}.
We use induction on $i$.
By the symmetry of $\left|\mathcal{M}_{2n+1,2i}\right|$ with
respect to $i$,
it suffices to consider the case  $i\leq (n+1)/2$.
Note that the theorem
 holds for $i=1$. Suppose that
$$\left|\mathcal{M}_{2n+1,2(i-1)}
\right|={n-1\choose i-2}{2n+1\choose n-1}.$$
By Theorems \ref{main}  \ref{syt}, we deduce that
\begin{align*}
\left|\mathcal{M}_{2n+1,2i}\right|&=\left|\mathcal{M}_{2n+1,2(i-1)}\right|\
+|\mathcal{T}_{2n+1,i}|\\[8pt]
&=\left({n-1\choose i-2}+\frac{n-2i+2}{i-1}{n-1\choose i-2}\right){2n+1\choose n-1}\\[8pt]
&={n-1\choose i-1}{2n+1\choose n-1},
\end{align*}
as desired. This completes the proof.
\qed

\vspace{.2cm} \noindent{\bf Acknowledgments.} We wish to thank
Sherry H.F. Yan for valuable comments. This work was
supported by  the 973 Project, the PCSIRT Project of the Ministry of
Education,  and the National
Science Foundation of China.

\end{document}